\documentstyle[12pt,a4,amscd,amssymb,stmaryrd,verbatim,graphics,]{amsart}
\newtheorem{thm}{Theorem}

\newtheorem{prop}{Poposition}
\newtheorem{ex}{Example}

\newtheorem{cor}{Corollary}
\newtheorem{defi}{Definition}

\begin{document}

\title{\bf The Lascoux, Leclerc and Thibon \\ algorithm and Soergel's 
tilting algorithm } 
\author{Steen Ryom-Hansen} 
\footnote{Supported in part by Programa Reticulados y Ecuaciones and by FONDECYT grant 1051024.} 
\address{Instituto de Matem\'atica y F\'isica \\
Univerdad de Talca \\
Chile \\ steen@@inst-mat.utalca.cl }
\begin{abstract}
{
We generalize Soergel's tilting algorithm to singular weights and 
deduce from this the validity of the Lascoux-Leclerc-Thibon conjecture on the connection between the 
canonical basis of the basic submodule of the Fock module and the representation theory of the Hecke-algebras at
root of unity.
}
\end{abstract}

\maketitle

\section {\bf Introduction }
\noindent

In this paper we show that the Lascoux-Leclerc-Thibon conjecture 
[LLT] on the connection between 
the canonical basis of the basic submodule of the Fock module and the 
representation theory of Hecke-algebras at a root of unity 
follows from the truth of Soergel's 
tilting algorithm. This result was independently obtained by Goodman and 
Wenzl [GW] and has also been proved by Leclerc-Thibon [LT].

\medskip

Our proof (which has existed in various versions since 1997) differs 
in several ways from the above proofs, first of all it relies 
notationally as wells as philosofically directly 
on the principle of graded representation
theory as exposed in the paper of Andersen, Jantzen 
and Soergel [AJS]. Indeed, as our first result   
we explain how the AJS-formalism 
naturally leads to an extension of Soergel's algorithm so as to be
able to deal with 
singular weights, i.e. weights lying on several reflecting
hyperplanes. This {\it singular combinatorics} is important from our
point of view, 
since the partitions appearing in the LLT-algorithm typically 
correspond to very singular weights. However, it should be noted that since the 
basic setting of [AJS] is that of Frobenius kernels, we cannot formally 
use the results of that paper. Indeed, we only check in type A that the singular combinatorics does not depend on 
the path of weights chosen and do so by comparing it with the LLT algorithm 
for large values of $l$, the order of the root of unity.

\medskip

We then go on to show that our singular combinatorics yields the correct tilting characters 
using the correctness of the original Soergel algorithm
together with some known properties of tilting modules and translation
functors.  
Finally, we show that the LLT-algorithm is a special case of our singular
combinatorics; this involves a detailed analysis of the correspondence
between partitions and weights. 

\medskip

I wish to thank  B. Leclerc and W. Soergel for useful discussions.


\section{\bf Preliminaries }
In this section we shall setup the notation needed. Let $ {\frak g}  $
be a finite dimensional semisimple Lie-algebra over the complex numbers 
and let $ U_q({\frak g} ) $ be 
the associated quantum group at an $ l$-th root of unity, see e.g. [A]
for the precise definition. 
The representation theory of $ U_q({\frak g} ) $ is 
labeled by the set of dominant weights $P^+ $ and the blocks 
correspond under this labeling to orbits in $ P^+ $ under the affine Weyl group $ W_l $.
Thus, for every weight $ \lambda \in P^+$ there is a standard module 
$\Delta(\lambda)$, a costandard module $\nabla(\lambda)$, a simple module 
$ L(\lambda) $ and a tilting module 
$Q(\lambda)$. See [A, AJS, S] for more details.

\medskip

We shall make extensive use of the following notation on alcove
geometry introduced in [AJS]. 
Let $ \Omega $ be a regular orbit of $ W_l $ in $ P^+ $, i.e. one
consisting of regular weights (lying on no walls), and let $ \Gamma $
be a singular orbit. Following [AJS], for $ \lambda \in \Omega $ 
we denote by $\lambda_{\Gamma} $ the unique element of $ \Gamma $ in
the closure of the alcove of $ \lambda $. 
Furthermore, we set 
$$ u(\lambda, \Gamma ) = | \{ H \in { \cal H } | \lambda_{\Gamma} \in
H, \lambda < H \} | $$ 
$$ o(\lambda, \Gamma ) = | \{ H \in { \cal H } | \lambda_{\Gamma} \in
H, \lambda > H \} | $$ 
where $ \cal H $ denotes the set of reflecting hyperplanes for $ W_l $.
Using this we can define $ \Omega_s $ by
$$ \Omega_s:= \{ \lambda \in \Omega \mid u(\lambda, \Gamma)=0 \} $$
Let us illustrate this on an A2 example:

\begin{ex}
\end{ex}
\begin{center}
\includegraphics{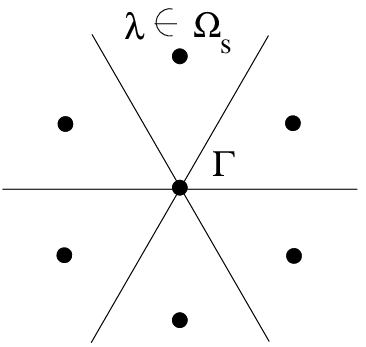}
\end{center}


\medskip
\noindent

Let us briefly recall Soergel's algorithm as well as the LLT-algorithm. Let 
$ \cal A $ be the set of alcoves, 
$ \cal A^+ $ the set of dominant alcoves. 
Then the Grothendieck group of $ \Omega $ can be identified with 
${\mathbb Z} [{\cal A^+}] $.
Soergel's algorithm produces for each $ A
\in \cal A^+ $ an ``indecomposable pattern'' $ P_A \in {\mathbb Z}[q] [{\cal A}^+]$, 
by which we mean an element $ P_A \in {\mathbb Z}[q] [{\cal A}^+]$ on the form 
$$ P_A(q) = A + \sum_{ B < A } \,  P_{AB}(q) $$
where $ P_{AB}(q) \in q {\mathbb Z}[q] $.
This pattern contains information about the character
of the tilting module $ Q(A) $ with highest weight $\lambda $ where $
\lambda \in A \cap \Omega $. In
formulas: 
$$ [Q(A),\Delta(B)]= P_{AB}(1) $$
where $ \Delta(B) $ 
is shorthand for $ \Delta(\mu) $, where $ \mu \in B \cap \Omega  $. 
The procedure for calculating $P_A$ is a recursion on $ \cal A $. 
It involves for each wall $s$ of the fundamental alcove
an operator $ \Theta_s $ on $ {\mathbb Z}[q] [{\cal A^+}] \, $ 
taking
the indecomposable pattern to a sum of indecomposable patterns (i.e. patterns with more
than one coefficient having a constant term). 
This operator is defined through the formula 
$$ 
\Theta_s A  = \left\{
\begin{array}{ll}
A + q (A s) & \mbox{if} \,\, As > A  \mbox{ and  } A s \in {\cal A}^+ \\
A + q^{-1}(A) & \mbox{if} \,\, A s < A \mbox{ and } As \in {\cal A}^+ \\
0 & \mbox{if} \,\, As \not\in {\cal A}^+ 
\end{array}
\right.
$$
and linearity; $ As $ is here the mirror alcove of $ A $ under the reflection given by $s$.
The following picture illustrates the first two cases of this action in a alcove geometry of an type A situation, where the reflection is
going upwards in the first case, downwards in the second case.

\begin{ex}
\end{ex}
\begin{center}
\includegraphics{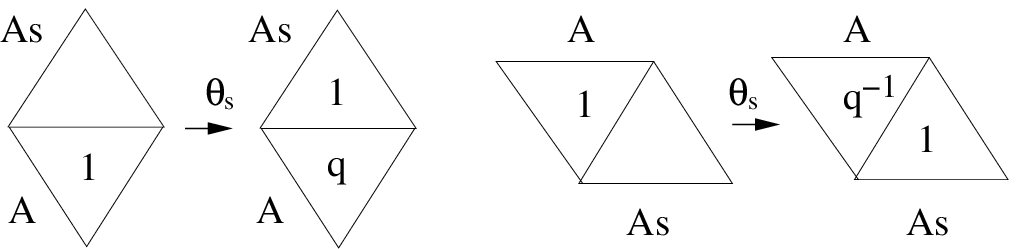}
\end{center}

One then 
subtracts inductively known indecomposable patterns to arrive at the new indecomposable pattern, whose top 
alcove is the only one with a coefficient involving a constant term.
We illustrate the algorithm on the following A3 examples,
and refer to [S1] for more details.

\begin{ex}
\end{ex}
\begin{center}
\includegraphics{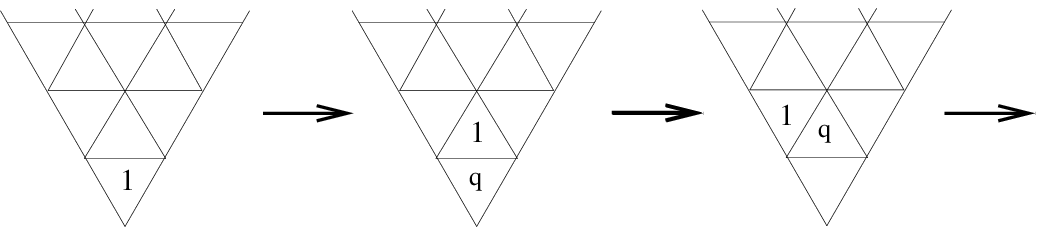}
\end{center}

\medskip

\begin{center}
\includegraphics{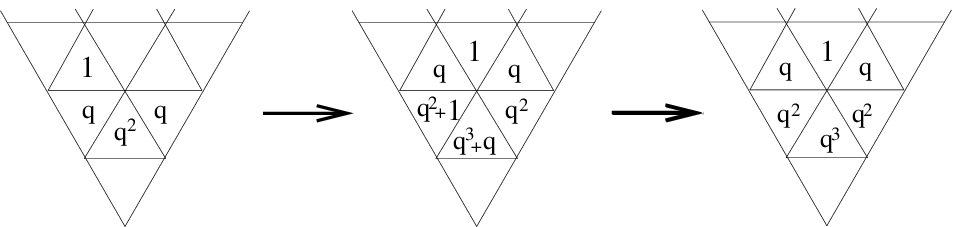}
\end{center}

\medskip

Following the terminology introduced in [S1], we shall denote the above algorithm a ``combinatorics'' for 
tilting modules.

\medskip
\medskip
We now briefly recall the LLT-algorithm. 
Let 
$ {\cal F}_q = \bigoplus_{\lambda \in { \rm Par}}
{\mathbb Q}(q ) | \lambda \rangle $ be the $q$-Fock space with basis 
parameterized by the set of all partitions Par. It can be made into
an integrable module for $U_q( \widehat{{\mathfrak sl}}_l)$ and thus
has a crystal basis. The LLT-algorithm calculates the global basis of
the basic submodule $ M$ of $ {\cal F}_q $, which is the one generated
by the empty partition. 
Let $ L $ be the ${\mathbb Z}[q]$-sublattice of ${\cal F}_q $
with basis $\{| \lambda \rangle \, | \, \lambda \in \mbox{Par} \} $.
The lower global basis element
$ G(\lambda ) $ of $ M $ is characterized by the following conditions
\begin{equation}\label{crystal}
\overline{G(\lambda)}  = G(\lambda), \quad G(\lambda) =
|\lambda\rangle 
\bmod qL 
\end{equation}
for $ \lambda \in \mbox{Par}_l $, i.e. an $l$-regular partition, 
where $ \overline{\, \cdot \, } $ is the involution of $ M $ given by 
$$ 
\begin{array}{ccccc}
\overline{\emptyset}=\emptyset, \, 
& \overline{f_i w}=
f_i \overline{ w } & \forall i  \, 
\mbox{ and } &
\overline{q}=q^{-1}, \, \, &
\end{array}
$$

\noindent
Let $ d_{\lambda \mu }(q) $ be defined by 
\begin{equation}\label{KL} G(\lambda) =\sum_{\mu  }
\, d_{\lambda \mu} (q) | \mu \rangle 
\end{equation}
Then $ d_{\lambda \mu}(q) 
\in {\mathbb Z}
[q] $, $ d_{\lambda,\mu}(q) = 0$ unless $ \lambda \trianglelefteq \mu $ and 
$ d_{\lambda,\lambda }(q)= 1 $.
Call an element $ w $ of $ M  $ selfdual if it satisfies
$ \overline{w} = w $.
A selfdual element $w $ can be written in the form 
$$ w = \sum_{\lambda} a_{\lambda}(q) G({\lambda}) $$
for some $ a_{\lambda}(q) \in {\mathbb Z}[q,q^{-1}] $ satisfying
$\overline{ a_{\lambda}(q)}= a_{\lambda}(q) $. The LLT-algorithm first 
constructs for each regular partition $ \lambda $ a selfdual element $ w_{\lambda}$
such that the coefficient of $ G(\lambda ) $ in $ w_{\lambda}$ is $1$ and
such that $ \mu <\lambda $ for all other occurring $ G(\mu ) $.
From this, $ G(\lambda) $ is obtained by linear 
algebra.

\medskip

LLT conjectured that for $ \lambda $ an $l$-regular partition
$$ d_{\lambda \mu}(1) = \left[ S(\mu), D(\lambda)  \right] $$ 
where $ S(\mu) $ and $ D(\lambda)$ are the Specht and 
the simple modules for the Hecke algebra of type A
specialized at an $ l$'th root of unity. This conjecture 
was first proved by Ariki [Ar] using the geometric approach 
to the crystal/canonical basis. The goal of this paper, however,
is to demonstrate that it also follows from Soergel's algorithm. 

\section{\bf Singular tilting modules }
 
We first need to generalize some results of Andersen on singular 
tilting modules. 



\medskip
Let $ \Omega $ be a regular $W_l$-orbit in $ P^+$ containing 
$ \lambda $ and let $ \Gamma $ be 
a singular orbit containing $ \mu $. Let $T_{\Gamma}^{\Omega}$ 
be the Jantzen translation functor from the $ \Gamma $-block
to the $\Omega $-block, see e.g. [A].
We then have the following proposition
\begin{prop}
Assume $ \lambda_{\Gamma} = \mu $ and  
$ \lambda \in \Omega_s $. Then 
$$ T_{\Gamma}^{\Omega} \, Q(\mu ) \cong Q(\lambda) $$
\end{prop}
\noindent {\it Proof}: One can copy the proof of Proposition 5.6 in Andersen's paper
[A]. In that paper $ \Gamma $ is assumed semiregular; however the
proof carries over to our situation.
\begin{flushright} $ \Box $ \end{flushright}

As a corollary, we obtain that the character of the singular tilting modules 
can be calculated from the regular ones:

\begin{cor} Let $ \mu, \overline{\mu} \in \Gamma $. Then 
$$[Q(\mu),\Delta(\overline{\mu}] = \frac{1}{N_{\Gamma}} 
\sum_{\overline{\lambda}:\overline{\lambda}_{\Gamma}= \overline{\mu}} \, 
[Q(\lambda), 
\Delta(\overline{
\lambda})] $$
where $ \lambda \in \Omega_s $ with $ \lambda_{\Gamma} = \mu $.
(As in [AJS] $ N_{\lambda } $ denotes the number of 
hyperplanes in $ \cal H $ such that $ \lambda \in H $).
\end{cor}
\noindent
{\bf Proof}: We know that $$ T^{\Gamma}_{\Omega} \, T_{\Gamma}^{\Omega}\,Q(\mu)
\cong Q(\mu)^{\oplus N_{\Gamma}}. $$ Hence we obtain the Corollary from 
the Theorem using 
the standard properties of translation functors.
\begin{flushright} $ \Box $ \end{flushright}

\section{\bf A combinatorics of graded translation functors}
We saw in the previous section that the singular tilting characters can be 
deduced from the regular ones. Now there may be no regular weights 
in the weight lattice, so we still insist on constructing a combinatorics 
of graded translation functors that works in the singular case. This
section is devoted to that task. We do it by assuming the existence of a 
formalism of graded translation functors for our $U_q({\frak g})$-representation theory 
having the same formal properties as the ones in [AJS] for the Frobenius 
kernels. We then show that this naturally leads to a combinatorics for singular 
tilting modules.

\medskip

Let thus $ \Omega_{\tau}= \Omega $, $ \Gamma_{\lambda}= \Gamma $ and 
$ \Pi_{\mu} = \Pi$ be the
orbits under $ W_l $ of $ \tau $, $ \lambda $ and $ \mu $. Assume furthermore 
$ \tau $ regular and $ \mu $ more singular than $ \lambda $, i.e.
$ W_{\lambda} \subseteq W_{\mu} $.

\medskip
\medskip
\begin{ex}
\end{ex}
\begin{center}
\includegraphics{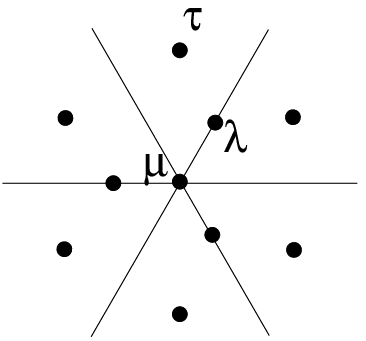}
\end{center}

\medskip

Let $ {\cal K}(\tau) $, $ {\cal K}(\lambda) $ and 
$ {\cal K}(\mu) $ be the corresponding $ \mathbb Z $-graded 
Grothendieck groups, i.e. $ {\cal K}(\tau) = {\mathbb Z }[q]
[\Omega]$ etc. As mentioned above they so far have been 
constructed only in the case of Frobenius kernels. 
According to the [AJS]-philosophy, they should come with a system of 
operators (graded translation functors)
$$ T^*_{\mu,\lambda}: {\cal K}(\mu) \rightarrow {\cal K}(\lambda), \,\,\,\,
\,\,\,\, T^*_{\lambda,\tau}:
{\cal K}(\lambda) \rightarrow {\cal K}(\tau), \,\,\,\,\,\,\,\,
T^*_{\mu,\tau}:
{\cal K}(\mu) \rightarrow {\cal K}(\tau) $$
$$ T_*^{\lambda,\mu}: {\cal K}(\lambda) \rightarrow {\cal K}(\mu), \,\,\,\,
\,\,\,\, T_*^{\tau,\lambda}:
{\cal K}(\tau) \rightarrow {\cal K}(\lambda), \,\,\,\,\,\,\,\,
T_*^{\tau,\mu}:
{\cal K}(\tau) \rightarrow {\cal K}(\mu) $$
as well as a system of $ T^! $ and $ T_! $ operators and a duality $ D $ 
relating the operators as follows
$$ D \circ T^*_{\mu,\lambda} \circ  D = T^!_{\mu,\lambda} $$ 
etc. Furthermore the duality should anticommute with the 
$ \mathbb Z $-shift in the categories, i.e. $ D \circ \langle 1 \rangle = 
\langle -1 \rangle \circ D $. 

\medskip

Now [AJS] page 253 suggests that $ T^*_{\mu,\tau} $ and $ T_*^{\tau,\mu} $ 
should satisfy the following rules:
$$ T_*^{\tau,\mu}\,\Delta(\tau) =  \Delta(\tau_{_{\Pi}})\,
\langle o(\tau, \Pi  ) \rangle $$
$$ T^*_{\mu,\tau}\,\Delta(\mu) = \sum_{\tau: \, \tau_{_{\Pi}}= \mu } \, \Delta(\tau)
\langle o(\tau, \Pi  ) \rangle $$
and similarly for $ T_{\lambda, \tau}^* $ and $ T_*^{\tau, \lambda} $. We
take this as our definition. 

\medskip
But then the transitivity forces us to define 
$ T_*^{\lambda, \mu } $ and $ T^*_{\mu, \lambda } $ by

$$  T_*^{\lambda, \mu } \, \Delta(\lambda) = \Delta( \lambda_{\Pi }) \, \langle o(\lambda, \Pi )
\rangle $$
$$ T^*_{ \mu, \lambda }\,\Delta(\mu) = \sum_{\lambda: \, 
\lambda_{_{\Pi}} = \mu } \, \Delta(\lambda) 
\langle o(\lambda, \Pi ) \rangle $$
i.e. the very same formulas as translation to and from the regular orbits. 
One should here notice that the expression $ \langle o(\lambda, \Pi)
\rangle $ makes sense for all weights and that 
$$ o(\tau, \Omega_{\lambda} ) + o(\tau_{_{\Gamma_{\lambda}}}, \Pi_{\mu})= 
o(\tau, \Pi_{\mu}) $$

Now $ T_*^{\lambda,\mu} =  T_!^{\lambda,\mu} $; 
hence  $ T_*^{\lambda,\mu} $ should preserve selfduality 
(i.e should commute with $ D $).
On the other hand we have that 
$ T^!_{\mu, \lambda} =  T^*_{\mu, \lambda}\, 
\langle -2(N_{\mu} - N_{\lambda})
\rangle $, so the operator that preserves selfduality should be 
$$ T^*_{\mu, \lambda} \, \langle N_{\lambda} - N_{\mu} \rangle $$
 
\medskip
Using the convention that $ \Delta( \lambda ) = 0 $ whenever $ \lambda 
\notin P^+ $ and that 
$$ o(\lambda,\Pi) + u(\lambda, \Pi ) = N_{\mu} - N_{\lambda} $$
\noindent
we arrive at the following first step for our combinatorics for tilting 
modules. We first assume that the character of $ Q(\nu) $ for 
$ \nu \in {\cal K}(\lambda) $ 
comes from an ``indecomposable pattern''
i.e. an element of $ {\mathbb Z}[q] {\cal K}(\lambda) $ on the form
$$  \nu+ \sum_{ \nu^{\prime} < \nu} \,  P_{\nu^{\prime},\nu}(q) \, \nu^{\prime} $$
with $ P_{\nu^{\prime},\nu}(q)  \in  q {\mathbb Z}[q] $, 
We then assume that there are operators akin to the 
$ \Theta_s $ of the Section 2. The above considerations lead 
us to choosing these as follows.

\medskip
\begin{defi} Singular combinatorics for tilting modules: Step 1:

\medskip

Let $ \lambda $ and $ \mu $ be as above, i.e. with $ \mu $ more singular than $ \lambda$ and
let $ {\cal K}(\mu)$, $ {\cal K}(\lambda)$ be the corresponding graded categories.
Then the graded translation functors $ \Theta^* $ and $ \Theta_* $ that take tilting modules to tilting modules, are: 

$$ \Theta_*:{\cal K}(\lambda) \rightarrow {\cal K}(\mu): \Delta(\lambda)
\mapsto \Delta(\lambda_{\Pi}) \, \langle o(\lambda,\Pi ) \rangle $$
$$ \Theta^*: {\cal K}(\mu) \rightarrow {\cal K}(\lambda):\,\,\,\,\,
\Delta(\mu) \mapsto  \sum_{\begin{array}{c} { \tiny
\lambda:  \lambda_{_{\Pi}} = \mu } \end{array}}  \, \Delta(\lambda) 
\langle -u(\lambda, \Pi ) \rangle $$
\end{defi}
\noindent

Let us illustrate this definition on an example

\begin{ex}
\end{ex}

\begin{center}
\includegraphics{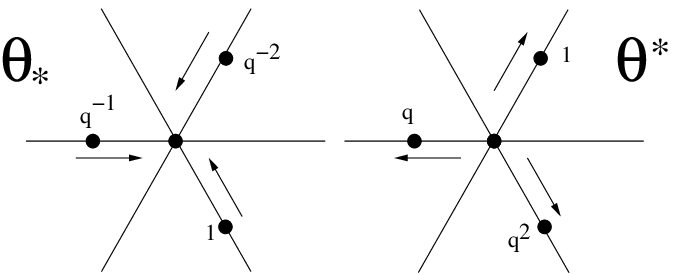}
\end{center}

\medskip
\noindent
We here used the convention that the shift $ \langle 1 \rangle $ in the 
graded category $ { \cal K }( \lambda ) $ corresponds to the multiplication
by $ q^{ -1} $.

\medskip
\medskip

The next step of our combinatorics is to explain how to obtain the
indecomposable pattern with highest weight $\mu $.

\medskip

This is, like in the regular case, an inductive procedure, starting with the weights $ \nu$ 
of the fundamental alcove, for which the pattern $ P_{\nu}(q) $ equals 
$ \nu $ itself. We then work ourselves upwards through 
the weight lattice with successive functors $ \Theta^* $ and $ \Theta_* $ always
trying to produce indecomposable´ patterns.

\medskip 

If $ P(\mu ) $ is an indecomposable pattern in $ {\cal K }(\mu) $ 
then it is clear from the definition that $ \Theta^* P(\mu) $ will 
remain indecomposable. 

\medskip
Now applying $ \Theta_* $ to an indecomposable pattern $ P(\lambda) $ in
$ {\cal K }(\lambda) $ will generally not produce an indecomposable pattern -- 
and $ \Theta_* P(\lambda) $ will generally not even have coefficients in 
$ {\mathbb Z } [q] $. We can therefore not just mimic Soergel's procedure
of subtracting inductively known patterns to arrive at something 
indecomposable.

\medskip

On the other hand the coefficient of the leading (maximal) weight $\mu $ will 
be ``1'' since $ \Theta_* $ does not lower the $ q$-power when going upwards. 
For each occurrence in the arising pattern of a $ q^{i} \nu $ with $ i $ negative or
zero, we then subtract $ \gamma(q) P_{\mu} (q) $ where 
$ \gamma(q) \in  {\mathbb Z } [q,q^{-1}] $ satisfies $ \gamma(q) = \gamma(q^{-1}) $. 
Repeating this eventually produces an indecomposable pattern with leading coefficient ``1''.



\medskip



\medskip
\medskip

Since there may be no weights inside the alcoves we need to generalize the 
previous considerations slightly to know what happens when translating
between arbitrary singular blocks.
Since translation functors should depend only on the alcove geometry, we 
may pretend that there exist regular weights. Then one finds for any 
$ \lambda, \lambda^{\prime} \in P^+ $ a $ \mu $ such that 
$$ W_{\lambda} \supseteq W_{\mu}, \,\,\, 
W_{\lambda^{\prime}} \supseteq W_{\mu}, \,\,\,
W_{\mu} \supseteq W_{\lambda} \cap W_{\lambda^{\prime}} $$

\medskip
Let us illustrate this situation on an example

\noindent
\begin{ex}
\end{ex}
\begin{center}
\includegraphics{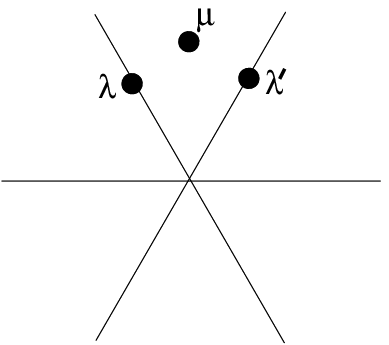}
\end{center}
\medskip

Now  
for the ordinary translation functors we have transitivity in that case:
$ T_{\lambda}^{\lambda^{\prime} } = T_{\mu}^{\lambda^{\prime}}
\circ T_{\lambda}^{\mu} $. We can therefore take the composite $ \Theta_{*} \circ \Theta^{*} $
as the graded version of $ T_{\lambda}^{\lambda^{\prime} } $.

\medskip 
\noindent
\begin{ex}
\end{ex}

\begin{center}
\includegraphics{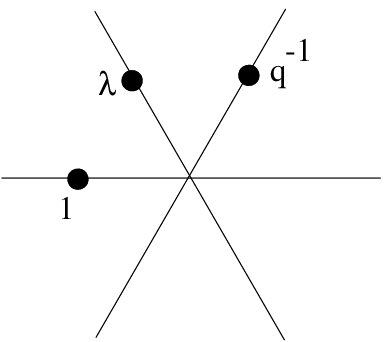}
\end{center}
\medskip

Finally, we kill all weights  which are not in the dominant Weyl chamber; this 
is analogous to Soergel's algorithm. Let us formulate all of this in one 
statement 

\medskip

\begin{defi} Singular combinatorics, step two.
Given $ \lambda \in P^+$. Let $ \nu_1,  \nu_2 , \ldots \nu_N  $ 
be the set of weights in $ P^+$ strictly less than $ \lambda $ in the usual 
order and assume inductively given tilting patterns 
$ P_{\nu_{i}}(q) \in {\mathbb Z}[q] {\cal A} $ for each $ \nu_i $. 
Let $ \nu  $ be a $ \nu_i$ in the closure of the alcove of $ \lambda $. 
Perform the relevant functor $ \Theta_* $ or $ \Theta^* $ or composite thereof on $ P_{\nu}(q) $
and subtract appropriate $ P_{\nu_{i}}(q) $'s as described to arrive at an 
indecomposable pattern. This is $ P_{\lambda}(q) $.

\end{defi} 

\section{\bf Comparing the combinatorics }

We now have a singular alcove combinatorics. It is clear that 
it gives Soergel's 
combinatorics if we only use semiregular orbits and translate through the 
walls. We must check that it always leads to the same answer, 
independently of the chosen path of weights. Once this has been established,
the algorithm will be correct, since we can choose a path 
$$ \nu_{1}, \nu_{2}, \nu_{3}, \ldots \, \nu_N, \nu_{N+1}, \ldots \nu_{N+K} $$
such that $ \nu_{1}, \nu_{2}, \nu_{3}, \ldots \, \nu_N $ are regular 
and semiregular while $ \nu_{N+1}, \ldots \nu_{N+K} $ have increasing
stabilizers. And for such a path, our algorithm yields the 
correct answer, by the correctness of Soergel's algorithm together with 
Corollary 1 and the construction of $\Theta_*$.



\medskip

We check this independency in type A only. The idea is to identify the graded translation
functors with the action of the $ f_i $'s on the Fock space; thus our 
combinatorics is really the combinatorics that calculates the global 
crystal basis. Since the global crystal basis is unique, the 
singular combinatorics will have no ambiguity either.

\medskip

Let us now therefore briefly review the correspondence between Young 
diagrams and weights in 
type $ {\rm A}_n $. 

\medskip

Let $ \lambda_i $ be the length of the $ i $'th line of the Young diagram
$ Y(\lambda) $. 
Then $Y(\lambda) $ is associated with the weight
$$ \lambda = ( \lambda_1 - \lambda_2, \lambda_2 - \lambda_3, \ldots, 
\lambda_{n-1} - \lambda_n, \lambda_n ) \in P^+ $$
due to the fact that the simple root $ \alpha_i $ in type $ A_n $
has coordinates 
$$ \alpha_i = (\,0, 0, \ldots, \stackrel{\,\,\downarrow i }{1}, \!\!\!\!
\!\!
\stackrel{\,\,\,\,\,\,\,\,\,\,\,\,\,\,\,\,\,\,\downarrow i+1 }{-1} , 0, 0 , 
\ldots \, ) $$

We conclude that $ \lambda + \rho $ lies on a wall corresponding to 
$ \alpha_i $ iff the last residues of the $ i $'th and the $ i+1 $'th 
rows are equal. This generalizes to other roots: if the last residues
of two rows are equal, the weight will lie on a wall. If one ends up 
on a wall by removing one node from the upper of two lines,
then the corresponding 
weight lies above the wall and so on. 

\medskip

\medskip
\noindent
\begin{ex}
\bf ($l=3$)
\end{ex}
\medskip
\medskip
\medskip
\medskip
\medskip
\medskip
\medskip
\medskip
\medskip
\medskip
\medskip
\medskip

\begin{center}
\begin{picture}(0,0)
\unitlength0.5cm

\multiput(-2,-1)(1,0){1}{\framebox(0.9,0.9){$0$}}
\multiput(-2,0)(1,0){1}{\framebox(0.9,0.9){$1$}}
\multiput(-2,1)(1,0){1}{\framebox(0.9,0.9){$2$}}
\multiput(-2,2)(1,0){1}{\framebox(0.9,0.9){$0$}}
\multiput(-2,3)(1,0){1}{\framebox(0.9,0.9){$1$}}
\multiput(-2,4)(1,0){1}{\framebox(0.9,0.9){$2$}}
\multiput(-2,5)(1,0){1}{\framebox(0.9,0.9){$0$}}

\multiput(-1,1)(1,0){1}{\framebox(0.9,0.9){$0$}}
\multiput(-1,2)(1,0){1}{\framebox(0.9,0.9){$1$}}
\multiput(-1,3)(1,0){1}{\framebox(0.9,0.9){$2$}}
\multiput(-1,4)(1,0){1}{\framebox(0.9,0.9){$0$}}
\multiput(-1,5)(1,0){1}{\framebox(0.9,0.9){$1$}}

\multiput(0,2)(1,0){1}{\framebox(0.9,0.9){$2$}}
\multiput(0,3)(1,0){1}{\framebox(0.9,0.9){$0$}}
\multiput(0,4)(1,0){1}{\framebox(0.9,0.9){$1$}}
\multiput(0,5)(1,0){1}{\framebox(0.9,0.9){$2$}}

\multiput(1,2)(1,0){1}{\framebox(0.9,0.9){$0$}}
\multiput(1,3)(1,0){1}{\framebox(0.9,0.9){$1$}}
\multiput(1,4)(1,0){1}{\framebox(0.9,0.9){$2$}}
\multiput(1,5)(1,0){1}{\framebox(0.9,0.9){$0$}}

\multiput(2,2)(1,0){1}{\framebox(0.9,0.9){$1$}}
\multiput(2,3)(1,0){1}{\framebox(0.9,0.9){$2$}}
\multiput(2,4)(1,0){1}{\framebox(0.9,0.9){$0$}}
\multiput(2,5)(1,0){1}{\framebox(0.9,0.9){$1$}}

\end{picture}
\end{center}

\medskip
\medskip
\medskip
\medskip

The $ 0$-residues give rise to three walls containing this $\lambda + \rho $. The removable $1$ node 
means that $ \lambda + \rho $ is positioned above two walls and below one wall coming from the $ 0$-nodes.

\medskip

Let us now focus on the $ n_0 $ rows of $ Y $ having as last node a 
$ 0 $-node. Let $ Y^{\prime} $ be the Young diagram obtained from $ Y $ 
by adding one node to one of the rows (\nolinebreak such a $ Y^{\prime} $ 
may not
exist). 

\medskip

We need to recall some facts on the modular representation theory
of $ GL_m(k) $.

\medskip
\medskip

Let $ \Delta(\lambda ) $ be the Weyl module given by $ \lambda
\vdash n $. It is a module for $ GL_m(k )$ for any $ m \geq n $.   
According to the branching rule ($ \alpha = \lambda, \, \, \beta = 1 $ in 
(2.30) of [J]) we have the following identity in the Grothendieck 
group:

$$ \Delta(\lambda) \, \otimes E = \, \sum_{\begin{array}{c} {\rm Young} (\lambda)  
\subseteq {\rm Young} (\mu) \\ \mid  {\rm Young} (\mu) \setminus 
{\rm Young} (\lambda) \mid \, = 1 \end{array} } \, 
\!\!\!\!\!\!\!\!\!\!\!\!\!\!\!\! \!\!\!\! \Delta( \mu ) $$ 

Here $ E = \Delta(1 ) $, i.e. the natural module for $ GL_m(k)
$. Using Donkin's version of the Nakayama conjecture [D] we obtain:
$$(*) \,\,\,\,\, pr_{\lambda^{\prime}} ( \Delta(\lambda) \, \otimes \, E ) = 
\sum_{\begin{array}{c} {\rm Young} (\lambda)  
\subseteq {\rm Young} (\mu) \\ \mid  {\rm Young} (\mu) \setminus 
{\rm Young} (\lambda) \mid = \{ 1{\rm -node } \}  \end{array} } 
\!\!\!\!\!\!\!\!\!\!\!\!\!\!\!\!\!\!\!\!\!\!\!\!\!\!\!\!\!\!\!\!
\Delta( \mu ) $$ 

Here $ pr_{\lambda^{\prime} } $ denotes projection onto the block of 
$ \Delta(\lambda^{\prime}) $; this is thus a formula for $ T_{\lambda}
^{\lambda^{\prime} } $ in the Grothendieck group.

\medskip

We are now going to calculate a graded version of this formula, in other
words we are going to apply the operators $ \Theta_* $ and $ \Theta^*$ of the
singular combinatorics of the previous section.

\medskip

Consider firstly the situation where the $ n_0 $ `$0$'-nodes are all removable 
and assume furthermore that $ n_1 = 0 $; i.e. no row has a a `$ 1 $'-node 
at the end. 

\begin{ex}
$ (l = 3) $
\end{ex}

\medskip
\medskip
\medskip
\medskip
\medskip
\medskip
\medskip
\medskip
\medskip

\begin{center}
\begin{picture}(0,0)
\unitlength0.5cm

\put(-10,1)


\multiput(-5,1)(1,0){1}{$\lambda = $}

\multiput(-2,-1)(1,0){1}{\framebox(0.9,0.9){$0$}}
\multiput(-2,0)(1,0){1}{\framebox(0.9,0.9){$1$}}
\multiput(-2,1)(1,0){1}{\framebox(0.9,0.9){$2$}}
\multiput(-2,2)(1,0){1}{\framebox(0.9,0.9){$3$}}
\multiput(-2,3)(1,0){1}{\framebox(0.9,0.9){$0$}}

\multiput(-1,0)(1,0){1}{\framebox(0.9,0.9){$2$}}
\multiput(-1,1)(1,0){1}{\framebox(0.9,0.9){$3$}}
\multiput(-1,2)(1,0){1}{\framebox(0.9,0.9){$0$}}
\multiput(-1,3)(1,0){1}{\framebox(0.9,0.9){$1$}}

\multiput(0,0)(1,0){1}{\framebox(0.9,0.9){$3$}}
\multiput(0,1)(1,0){1}{\framebox(0.9,0.9){$0$}}
\multiput(0,2)(1,0){1}{\framebox(0.9,0.9){$1$}}
\multiput(0,3)(1,0){1}{\framebox(0.9,0.9){$2$}}

\multiput(1,0)(1,0){1}{\framebox(0.9,0.9){$0$}}
\multiput(1,1)(1,0){1}{\framebox(0.9,0.9){$1$}}
\multiput(1,2)(1,0){1}{\framebox(0.9,0.9){$2$}}
\multiput(1,3)(1,0){1}{\framebox(0.9,0.9){$3$}}

\multiput(2,0)(1,0){1}{\framebox(0.9,0.9){$1$}}
\multiput(2,1)(1,0){1}{\framebox(0.9,0.9){$2$}}
\multiput(2,2)(1,0){1}{\framebox(0.9,0.9){$3$}}
\multiput(2,3)(1,0){1}{\framebox(0.9,0.9){$0$}}

\multiput(3,0)(1,0){1}{\framebox(0.9,0.9){$2$}}
\multiput(3,1)(1,0){1}{\framebox(0.9,0.9){$3$}}
\multiput(3,2)(1,0){1}{\framebox(0.9,0.9){$0$}}
\multiput(3,3)(1,0){1}{\framebox(0.9,0.9){$1$}}

\multiput(4,1)(1,0){1}{\framebox(0.9,0.9){$0$}}
\multiput(4,2)(1,0){1}{\framebox(0.9,0.9){$1$}}
\multiput(4,3)(1,0){1}{\framebox(0.9,0.9){$2$}}

\multiput(5,2)(1,0){1}{\framebox(0.9,0.9){$2$}}
\multiput(5,3)(1,0){1}{\framebox(0.9,0.9){$3$}}

\multiput(6,3)(1,0){1}{\framebox(0.9,0.9){$0$}}

\multiput(6,3)(1,0){1}{\framebox(0.9,0.9){$0$}}

\end{picture}
\end{center}

\medskip
\medskip
\medskip
\medskip

The $ n_0 $ `$ 0 $'-nodes give rise to $ \left( \begin{array}{c} n_0 \\ 2 \end{array} 
\right) $ hyperplanes through $ \lambda $. The other residues also 
give rise to hyperplanes through $ \lambda $; on the other hand, 
the components
of $ T_{\lambda}^{\lambda^{\prime} } \, \Delta(\lambda) $ all stay fixed 
with respect to these
other hyperplanes since we are adding only `$1$'-nodes in $(*)$.

\medskip
 
The components
of $ T_{\lambda}^{\lambda^{\prime} } \, \Delta(\lambda) $ lie on fewer 
hyperplanes than $ \lambda  $ since $ n_1 = 0 $ and thus we are in position to 
apply $ \Theta_{\lambda, \lambda^{\prime} }^* $ from our singular combinatorics.

\medskip

We must for each component $ \Delta(\mu) $ of 
$ T_{\lambda}^{\lambda^{\prime} } \, \Delta(\lambda) $ calculate the number 
$$ u(\mu, \Gamma_{\lambda} ) = \# \{ \, \, {\rm hyperplanes \,\,\,} H \,\,\, 
{\rm through \,\,\, } \lambda \,\,\, {\rm with } \,\,\, \mu  < H \} $$
With the above $ \lambda $, one of the occurring $ \mu $'s will be 

\medskip
\medskip
\medskip
\medskip
\medskip
\medskip
\medskip
\medskip
\medskip
\medskip
\medskip

\begin{center}
\begin{picture}(0,0)
\unitlength0.5cm

\multiput(-12,-1)(1,0){1}{\framebox(0.9,0.9){$0$}}
\multiput(-12,0)(1,0){1}{\framebox(0.9,0.9){$1$}}
\multiput(-12,1)(1,0){1}{\framebox(0.9,0.9){$2$}}
\multiput(-12,2)(1,0){1}{\framebox(0.9,0.9){$3$}}
\multiput(-12,3)(1,0){1}{\framebox(0.9,0.9){$0$}}

\multiput(-11,0)(1,0){1}{\framebox(0.9,0.9){$2$}}
\multiput(-11,1)(1,0){1}{\framebox(0.9,0.9){$3$}}
\multiput(-11,2)(1,0){1}{\framebox(0.9,0.9){$0$}}
\multiput(-11,3)(1,0){1}{\framebox(0.9,0.9){$1$}}

\multiput(-10,0)(1,0){1}{\framebox(0.9,0.9){$3$}}
\multiput(-10,1)(1,0){1}{\framebox(0.9,0.9){$0$}}
\multiput(-10,2)(1,0){1}{\framebox(0.9,0.9){$1$}}
\multiput(-10,3)(1,0){1}{\framebox(0.9,0.9){$2$}}

\multiput(-9,0)(1,0){1}{\framebox(0.9,0.9){$0$}}
\multiput(-9,1)(1,0){1}{\framebox(0.9,0.9){$1$}}
\multiput(-9,2)(1,0){1}{\framebox(0.9,0.9){$2$}}
\multiput(-9,3)(1,0){1}{\framebox(0.9,0.9){$3$}}

\multiput(-8,0)(1,0){1}{\framebox(0.9,0.9){$1$}}
\multiput(-8,1)(1,0){1}{\framebox(0.9,0.9){$2$}}
\multiput(-8,2)(1,0){1}{\framebox(0.9,0.9){$3$}}
\multiput(-8,3)(1,0){1}{\framebox(0.9,0.9){$0$}}

\multiput(-7,0)(1,0){1}{\framebox(0.9,0.9){$2$}}
\multiput(-7,1)(1,0){1}{\framebox(0.9,0.9){$3$}}
\multiput(-7,2)(1,0){1}{\framebox(0.9,0.9){$0$}}
\multiput(-7,3)(1,0){1}{\framebox(0.9,0.9){$1$}}

\multiput(-6,1)(1,0){1}{\framebox(0.9,0.9){$0$}}
\multiput(-6,2)(1,0){1}{\framebox(0.9,0.9){$1$}}
\multiput(-6,3)(1,0){1}{\framebox(0.9,0.9){$2$}}

\multiput(-5,2)(1,0){1}{\framebox(0.9,0.9){$2$}}
\multiput(-5,3)(1,0){1}{\framebox(0.9,0.9){$3$}}

\multiput(-4,3)(1,0){1}{\framebox(0.9,0.9){$0$}}

\put(-9,-2){$ \lambda $ }

\multiput(3,-1)(1,0){1}{\framebox(0.9,0.9){$0$}}
\multiput(3,0)(1,0){1}{\framebox(0.9,0.9){$1$}}
\multiput(3,1)(1,0){1}{\framebox(0.9,0.9){$2$}}
\multiput(3,2)(1,0){1}{\framebox(0.9,0.9){$3$}}
\multiput(3,3)(1,0){1}{\framebox(0.9,0.9){$0$}}

\multiput(4,-1)(1,0){1}{\framebox(0.9,0.9){$1$}}
\multiput(4,0)(1,0){1}{\framebox(0.9,0.9){$2$}}
\multiput(4,1)(1,0){1}{\framebox(0.9,0.9){$3$}}
\multiput(4,2)(1,0){1}{\framebox(0.9,0.9){$0$}}
\multiput(4,3)(1,0){1}{\framebox(0.9,0.9){$1$}}

\multiput(5,0)(1,0){1}{\framebox(0.9,0.9){$3$}}
\multiput(5,1)(1,0){1}{\framebox(0.9,0.9){$0$}}
\multiput(5,2)(1,0){1}{\framebox(0.9,0.9){$1$}}
\multiput(5,3)(1,0){1}{\framebox(0.9,0.9){$2$}}

\multiput(6,0)(1,0){1}{\framebox(0.9,0.9){$0$}}
\multiput(6,1)(1,0){1}{\framebox(0.9,0.9){$1$}}
\multiput(6,2)(1,0){1}{\framebox(0.9,0.9){$2$}}
\multiput(6,3)(1,0){1}{\framebox(0.9,0.9){$3$}}

\multiput(7,0)(1,0){1}{\framebox(0.9,0.9){$1$}}
\multiput(7,1)(1,0){1}{\framebox(0.9,0.9){$2$}}
\multiput(7,2)(1,0){1}{\framebox(0.9,0.9){$3$}}
\multiput(7,3)(1,0){1}{\framebox(0.9,0.9){$0$}}

\multiput(8,0)(1,0){1}{\framebox(0.9,0.9){$2$}}
\multiput(8,1)(1,0){1}{\framebox(0.9,0.9){$3$}}
\multiput(8,2)(1,0){1}{\framebox(0.9,0.9){$0$}}
\multiput(8,3)(1,0){1}{\framebox(0.9,0.9){$1$}}

\multiput(9,1)(1,0){1}{\framebox(0.9,0.9){$0$}}
\multiput(9,2)(1,0){1}{\framebox(0.9,0.9){$1$}}
\multiput(9,3)(1,0){1}{\framebox(0.9,0.9){$2$}}

\multiput(10,2)(1,0){1}{\framebox(0.9,0.9){$2$}}
\multiput(10,3)(1,0){1}{\framebox(0.9,0.9){$3$}}

\multiput(11,3)(1,0){1}{\framebox(0.9,0.9){$0$}}

\put(6,-2){$ \mu $ }

\end{picture}
\end{center}

\medskip
\medskip
\medskip
\medskip
\medskip

One sees in this example (which is easily generalized to all $ \mu $) that 
$$ u(\mu, \Gamma_{\lambda} ) = \# \{ 0\mbox{-nodes of } {\lambda}
\mbox{ above } 
\gamma \}$$
where $ \gamma = \mu \setminus \lambda  $.
Thus $\Delta( \mu)  $ occurs in the graded translation from $ \lambda $ to $ \lambda^
{\prime} $ with a shift of order:
$$ \begin{array}{l}
 - u(\mu, \Gamma_{\lambda} ) =  -  \# \{ 0\mbox{-nodes of } {\lambda}
 \mbox{ above }
\gamma  \}  \\   = 
-\# \{\mbox{ indent } 1 \mbox{-nodes of } { \lambda} \mbox{ above }
\gamma  \} +
\# \{\mbox{ removable } 1 \mbox{-nodes of } { \lambda } \mbox{ above }\gamma \}  $$
\end{array} $$

\medskip

We now consider the slightly more general situation where some of the `$ 0 $'-nodes
of $ \lambda $ are allowed to be non-removable; on the other hand we will 
require that there be no removable `$ 1 $'-nodes in $ \lambda $. One then 
gets that 
$$ n_1 = \#\{ \mbox{ non-removable } 0\mbox{-nodes of } \lambda \} $$

\begin{ex}
$ ( l = 3) $ 
\end{ex}

\medskip
\medskip
\medskip
\medskip
\medskip
\medskip
\medskip
\medskip
\medskip
\medskip
\medskip

\begin{center}
\begin{picture}(0,0)
\unitlength0.5cm

\multiput(-8,0)(1,0){1}{\framebox(0.9,0.9){$2$}}
\multiput(-8,1)(1,0){1}{\framebox(0.9,0.9){$0$}}
\multiput(-8,2)(1,0){1}{\framebox(0.9,0.9){$1$}}
\multiput(-8,3)(1,0){1}{\framebox(0.9,0.9){$2$}}
\multiput(-8,4)(1,0){1}{\framebox(0.9,0.9){$0$}}

\multiput(-7,0)(1,0){1}{\framebox(0.9,0.9){$0$}}
\multiput(-7,1)(1,0){1}{\framebox(0.9,0.9){$1$}}
\multiput(-7,2)(1,0){1}{\framebox(0.9,0.9){$2$}}
\multiput(-7,3)(1,0){1}{\framebox(0.9,0.9){$0$}}
\multiput(-7,4)(1,0){1}{\framebox(0.9,0.9){$1$}}

\multiput(-6,1)(1,0){1}{\framebox(0.9,0.9){$2$}}
\multiput(-6,2)(1,0){1}{\framebox(0.9,0.9){$0$}}
\multiput(-6,3)(1,0){1}{\framebox(0.9,0.9){$1$}}
\multiput(-6,4)(1,0){1}{\framebox(0.9,0.9){$2$}}

\multiput(-5,4)(1,0){1}{\framebox(0.9,0.9){$0$}}

\put(-10,2){$ \lambda =$ }

\multiput(3,0)(1,0){1}{\framebox(0.9,0.9){$2$}}
\multiput(3,1)(1,0){1}{\framebox(0.9,0.9){$0$}}
\multiput(3,2)(1,0){1}{\framebox(0.9,0.9){$1$}}
\multiput(3,3)(1,0){1}{\framebox(0.9,0.9){$2$}}
\multiput(3,4)(1,0){1}{\framebox(0.9,0.9){$0$}}

\multiput(4,0)(1,0){1}{\framebox(0.9,0.9){$0$}}
\multiput(4,1)(1,0){1}{\framebox(0.9,0.9){$1$}}
\multiput(4,2)(1,0){1}{\framebox(0.9,0.9){$2$}}
\multiput(4,3)(1,0){1}{\framebox(0.9,0.9){$0$}}
\multiput(4,4)(1,0){1}{\framebox(0.9,0.9){$1$}}

\multiput(5,0)(1,0){1}{\framebox(0.9,0.9){$1$}}
\multiput(5,1)(1,0){1}{\framebox(0.9,0.9){$2$}}
\multiput(5,2)(1,0){1}{\framebox(0.9,0.9){$0$}}
\multiput(5,3)(1,0){1}{\framebox(0.9,0.9){$1$}}
\multiput(5,4)(1,0){1}{\framebox(0.9,0.9){$2$}}

\multiput(6,4)(1,0){1}{\framebox(0.9,0.9){$0$}}

\put(7,2){$= \mu $ }

\end{picture}
\end{center}

\medskip

The $ n_0 $ `$ 0 $'-nodes still give rise to 
$ \left( \begin{array}{c} n_0 \\ 2 \end{array} \right) $ walls passing 
through $ \lambda $, but in this situation $ \Theta^*_{\lambda, \lambda^{\prime}} $
is not a graded version of $ T_{\lambda}^{\lambda^{\prime} } $, since 
$ \lambda^{\prime} $ lies on walls not passing through $ \lambda $. 

\medskip
The correct graded version of $ T_{\lambda}^{\lambda^{\prime} } $ will 
be $ \Theta_*^{\tau,\lambda^{\prime}}  \circ  \Theta_{\lambda, \tau }^* $,
where $ \tau \in X(T)\otimes {\mathbb R } $ is more regular than each of the 
weights $ \lambda $ and $ \lambda^{\prime} $. In our example $ \tau $ can
be visualized by the following ``diagram''. 

\medskip
\medskip

\medskip
\medskip
\medskip
\medskip
\medskip
\medskip
\medskip
\medskip 

\begin{center}
\begin{picture}(0,0)
\unitlength0.5cm

\multiput(-3,0)(1,0){1}{\framebox(0.9,0.9){$2$}}
\multiput(-3,1)(1,0){1}{\framebox(0.9,0.9){$0$}}
\multiput(-3,2)(1,0){1}{\framebox(0.9,0.9){$1$}}
\multiput(-3,3)(1,0){1}{\framebox(0.9,0.9){$2$}}
\multiput(-3,4)(1,0){1}{\framebox(0.9,0.9){$0$}}

\multiput(-2,0)(1,0){1}{\framebox(0.9,0.9){$0$}}
\multiput(-2,1)(1,0){1}{\framebox(0.9,0.9){$1$}}
\multiput(-2,2)(1,0){1}{\framebox(0.9,0.9){$2$}}
\multiput(-2,3)(1,0){1}{\framebox(0.9,0.9){$0$}}
\multiput(-2,4)(1,0){1}{\framebox(0.9,0.9){$1$}}

\multiput(-1,1)(1,0){1}{\framebox(0.9,0.9){$2$}}
\multiput(-1,2)(1,0){1}{\framebox(0.9,0.9){$0$}}
\multiput(-1,3)(1,0){1}{\framebox(0.9,0.9){$1$}}
\multiput(-1,4)(1,0){1}{\framebox(0.9,0.9){$2$}}

\multiput(-1,0)(1,0){1}{\framebox(0.4,0.9)}

\multiput(0,4)(1,0){1}{\framebox(0.9,0.9){$0$}}

\put(-5,2){$ \tau = $ }

\end{picture}
\end{center}

\medskip
\noindent
where we have added ``half a node''.

\medskip

We then have $$ u(\tau ,\Gamma_{\lambda} ) = \#\{ 0\mbox{-nodes of } 
\lambda \,\,\, \mbox{above } \mu \setminus \lambda \} $$
Furthermore we have that 
$$ o(\tau ,\Gamma_{\mu} ) = \#\{ 1\mbox{-nodes of } \lambda \,\,\, 
\mbox{above } \mu \setminus \lambda \} $$
So $ \Delta(\mu) $ appears in the graded translation functor from $ \lambda $ 
to $ \lambda^{\prime} $ with a shift of size 
$$ -u(\tau, \Gamma_{\lambda}) + o(\tau, \Gamma_{\mu}) = 
- \#\{ { \rm indent \,\,\,} 1 \mbox{-nodes of } \lambda \mbox{ above } \mu \setminus \lambda \} $$
$$ = - \#\{ { \rm indent \,\,\,} 1 \mbox{-nodes of } \lambda \mbox{ above } \mu \setminus
\lambda \} +  \#\{ { \rm removable \,\,\,} 1 \mbox{-nodes of } \lambda
\mbox{ above } \mu \setminus \lambda \} $$

\medskip
\medskip
\medskip

Consider finally the general situation in which we allow removable 
`$ 1 $'-nodes.

\medskip 

Our graded version of $ T_{\lambda}^{\lambda^{\prime} } $ will then behave 
like $ \oplus \, \Theta_*^{\tau, \lambda^{\prime} } \circ \Theta_{\lambda, \tau}^*  $
where $ \tau \in X(T) \otimes {\mathbb R } $ is chosen as before. Each 
$ \Delta(\mu ) $ appearing in the graded translation corresponds to a $ \tau $ 
and the same calculation as before gives a shift of size 

$$ - \# \{\,\,\, { \rm indent } \,\,\, 1 \mbox{-nodes of } \lambda
\mbox{ above } \mu \setminus 
\lambda \,\,\,\} $$
$$ +\#\{\,\,\, { \rm removable }\,\,\, 1 \mbox{-nodes of } \lambda
\mbox{ above } \mu \setminus 
\lambda \,\,\, \} $$

\medskip 

Hence we get in all cases exactly the number $ -N_1^r( \lambda, \mu ) $ of LLT. We 
can of 
course repeat this argument for the other residues and we conclude 
that the $ f_i $ operators on the Fock space are really those graded 
translation functors in our singular setup that are summands of the tensor product
with the natural module. Let us formulate this as a theorem:

\begin{thm}
Let $ f_i $ be one of the standard generators of $  U_q( \widehat{ {\mathfrak
  sl }_l } ) $ and let $ f_i \lambda =  \sum c_{\mu}(q) \mu $ in the action
of $ f_i $ on $ {\cal F}_q $. Then 
$ \sum c_{\mu}(q) \mu $ equals $ \Theta  \Delta(\lambda) $ in the singular
combinatorics of the previous section, where $ \Theta $ is the operator of the singular
combinatorics corresponding to $ \lambda $ and $ \mu $, with $ \mu $
being obtained from $ \lambda $ by adding an $i$-node.

\end{thm}

\medskip 

Now recall the selfdual element $ w_{\lambda} $ of $ {\cal F}_q $, which 
is the first step of the LLT-algorithm. It is 
on the form
$$ w_{\lambda} = f_{i_1}^{(n_1)} f_{i_2}^{(n_2)} \ldots f_{i_k}^{(n_k)} \emptyset $$
for some $ i_k $ and $ n_k $ where $ f_i^{(n)} = ([n]_q!)^{-1} \, f_i^n  $ is the usual
divided power notation.



\medskip 
We therefore need to check that also the action of the higher 
divided powers  
$ f_i^{(n)} = ([n]_q!)^{-1} \, f_i^n  $ can be described in the
singular combinatorics. This is an argument close to the above. Let us
start with an example.

\begin{ex}
$(l=3)$
\end{ex}

\medskip
\medskip
\medskip
\medskip
\medskip
\medskip
\medskip
\medskip

\begin{center}
\begin{picture}(0,0)
\unitlength0.5cm

\multiput(-4,1)(1,0){1}{\framebox(0.9,0.9){$2$}}
\multiput(-4,2)(1,0){1}{\framebox(0.9,0.9){$3$}}
\multiput(-4,3)(1,0){1}{\framebox(0.9,0.9){$0$}}

\multiput(-3,1)(1,0){1}{\framebox(0.9,0.9){$3$}}
\multiput(-3,2)(1,0){1}{\framebox(0.9,0.9){$0$}}
\multiput(-3,3)(1,0){1}{\framebox(0.9,0.9){$1$}}

\multiput(-2,1)(1,0){1}{\framebox(0.9,0.9){$0$}}
\multiput(-2,2)(1,0){1}{\framebox(0.9,0.9){$1$}}
\multiput(-2,3)(1,0){1}{\framebox(0.9,0.9){$2$}}

\multiput(-1,2)(1,0){1}{\framebox(0.9,0.9){$2$}}
\multiput(-1,3)(1,0){1}{\framebox(0.9,0.9){$3$}}

\multiput(0,2)(1,0){1}{\framebox(0.9,0.9){$3$}}
\multiput(0,3)(1,0){1}{\framebox(0.9,0.9){$0$}}

\multiput(1,2)(1,0){1}{\framebox(0.9,0.9){$0$}}
\multiput(1,3)(1,0){1}{\framebox(0.9,0.9){$1$}}

\multiput(2,3)(1,0){1}{\framebox(0.9,0.9){$2$}}

\multiput(3,3)(1,0){1}{\framebox(0.9,0.9){$3$}}

\multiput(4,3)(1,0){1}{\framebox(0.9,0.9){$0$}}

\put(-6,2){$ \tau \,\,=  $ }

\end{picture}
\end{center}

\medskip
\medskip

It corresponds to a Steinberg weight in the $A_2 $ situation. 
Translating to the weight given by the diagram $ \sigma $:

\medskip
\medskip
\medskip
\medskip
\medskip
\medskip
\medskip

\begin{center}
\begin{picture}(0,0)
\unitlength0.5cm

\multiput(-4,1)(1,0){1}{\framebox(0.9,0.9){$2$}}
\multiput(-4,2)(1,0){1}{\framebox(0.9,0.9){$3$}}
\multiput(-4,3)(1,0){1}{\framebox(0.9,0.9){$0$}}

\multiput(-3,1)(1,0){1}{\framebox(0.9,0.9){$3$}}
\multiput(-3,2)(1,0){1}{\framebox(0.9,0.9){$0$}}
\multiput(-3,3)(1,0){1}{\framebox(0.9,0.9){$1$}}

\multiput(-2,1)(1,0){1}{\framebox(0.9,0.9){$0$}}
\multiput(-2,2)(1,0){1}{\framebox(0.9,0.9){$1$}}
\multiput(-2,3)(1,0){1}{\framebox(0.9,0.9){$2$}}

\multiput(-1,2)(1,0){1}{\framebox(0.9,0.9){$2$}}
\multiput(-1,3)(1,0){1}{\framebox(0.9,0.9){$3$}}

\multiput(0,2)(1,0){1}{\framebox(0.9,0.9){$3$}}
\multiput(0,3)(1,0){1}{\framebox(0.9,0.9){$0$}}

\multiput(1,2)(1,0){1}{\framebox(0.9,0.9){$0$}}
\multiput(1,3)(1,0){1}{\framebox(0.9,0.9){$1$}}

\multiput(2,2)(1,0){1}{\framebox(0.9,0.9){$1$}}
\multiput(2,3)(1,0){1}{\framebox(0.9,0.9){$2$}}

\multiput(3,3)(1,0){1}{\framebox(0.9,0.9){$3$}}

\multiput(4,3)(1,0){1}{\framebox(0.9,0.9){$0$}}

\multiput(5,3)(1,0){1}{\framebox(0.9,0.9){$1$}}

\put(-6,2){$ \sigma \,\,=  $ }

\end{picture}
\end{center}
can be described in our combinatorics by the following picture

\medskip
\medskip

\begin{center}
\includegraphics{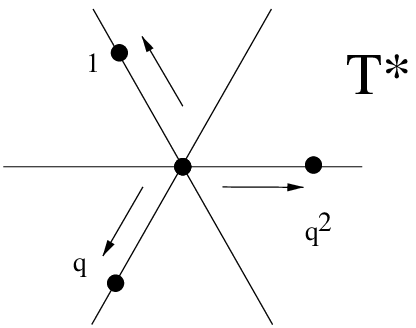}
\end{center}

\medskip
\medskip
\medskip

This is performed in the LLT-algorithm 
as a two step operation, adding one node in each step. In the 
alcove geometry we get the following picture
\medskip
\medskip

\begin{center}
\includegraphics{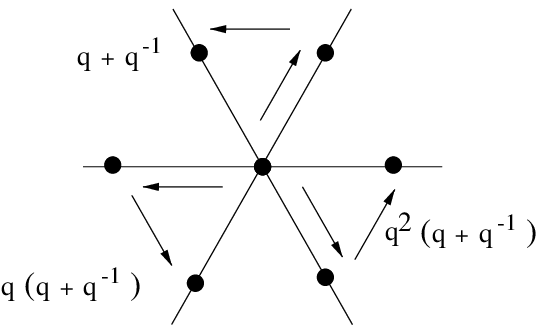}
\end{center}

\medskip
\medskip


\noindent
and so $ \Theta^*$ corresponds to 
$ f_1^{(2)} = \frac{1}{[2]} f_1^{2}$. 

\medskip

The example generalizes to all $ n$:

\begin{thm}
The action of $f_i^{(n)} $ on $ {\cal F}_q $ 
corresponds to a $ \Theta $ of the singular combinatorics. 
\end{thm}
{\bf Proof}: 
We may assume that $ i \equiv 1 \mod l $. 
Assume first that the residues of the end nodes of
$ \lambda $ are all different from $ 1 $. Then the end 
nodes of residue $ 0 $ all give rise to addable $1$-nodes.
Let us  
consider $f_1^{(n)} \lambda $. Let $ I = \{ i_{l_1}, i_{l_2}, \ldots ,  i_{l_K} \} $
be the lines of $ \lambda $ having $0$ as end residue and let $ \lambda^J $ 
for $ J \subset \{ l_1, \ldots , l_K \} $ with $|J|=n$ be the partition obtained from $\lambda $ by adding one node to the lines 
of $ i_{l_j}  $ for $ j \in J  $. The pairs $ (a,b) $ of line numbers
where $ a \in I $ such that $ a < j \,\,  \forall j \in J $ and $ b \in J $
correspond to the hyperplanes through $ \lambda $ that lie above 
$ \lambda^J $, i.e. those contributing to $ -u(\lambda^J, \Pi ) $ of
Definition 1 where $ \lambda \in \Pi $. But in the notation of Lemma 6.2 of [LLT] there are 
exactly $ N(id) + \left( \begin{array}{c} n \\ 2 \end{array} \right)
$ of these (where $ n = k_s $ in {\it loc. cit.}). Since this is 
also the coefficient of $ \lambda^J $ in $f_1^{(n)} \lambda $
the proof is finished in that case.  

\medskip

If there are $1$-nodes occurring at the end of some lines, the 
situation is slightly more complicated since the relevant operator
in the combinatorics is a composite of $ \Theta^* $ and $ \Theta_* $.
Let $J$ be as before.
Let us first assume that these $1$-nodes are all removable and let 
us denote by $K$ the corresponding line numbers.
Then the pairs $ (a,b) $ of line numbers
where $ a \in K $ with $ a < j \,\,  \forall j \in J $ 
and $b \in J $ gives the number of hyperplanes through $\lambda $ lying 
below $ \lambda^J$, i.e. contributing to $ o(\lambda^J, \Pi ) $.
By Definition 1 we should subtract this 
number from $ -u(\lambda^J, \Pi ) $
and thus find once more 
the exact correspondence with the formula for $f_1^{(n)} \lambda $.

\medskip
Actually, this argument also holds in the case where
some of these $1$-nodes are non-removable, by our definition 
of $\Theta $ via ``half''-nodes and we are done.

\begin{flushright} $ \Box $ \end{flushright}

We can now prove that the singular combinatorics is well defined:

\begin{thm}
Let the root system be of type A. Then
the singular combinatorics is well defined, i.e. does not depend on the
path of weights. 
\end{thm}
\noindent {\it Proof}: 
Let $ \mu^1, \mu^1, \ldots , \mu^k, \ldots,  \mu^N $ 
be a path of weights in $ P^+$, starting with $ \mu^1 $ in the 
fundamental alcove $ C $ 
and finishing with $ \mu^N = \lambda $. We view the weights as
partitions and show that each of 
the operators $ \Theta_{\mu^i, \mu^{i+1}} $ of 
the singular combinatorics applied to $ \mu^i $ can be identified with the action of 
$ f_{i_1}^{(n_1)} f_{i_2}^{(n_2)} \cdots f_{i_p}^{(n_p)} $ on $ \nu^i \in {\cal F}_q $ for
some choice of $ f_{i_j}  $ and $ n_{i_j}$. But then
the uniqueness of the crystal 
basis of $ M \subset {\cal F}_q $ shows
that the crystal basis is independent of the choice of path.

\medskip

Let $ T \subset G = Sl_m $ be a maximal torus so that $\mu^k $ defines a
$T$-module. We associate a partition $ (\mu^i_1, \mu^i_2, \ldots,
\mu^i_m) $
with each of the $ \mu^i $ by the rule $ \langle \mu^i, \alpha_1
\rangle =   \mu^i_1- \mu^i_2 $ etc. This partition is unique up to adding the first columns
a number of times.

\medskip 
 
As usual, we let $ \rho $ denote the half sum of the positive roots, the corresponding
partition being $  ( m, m-1, \ldots, 1 ) $. Recall that $ l $ is the
order of the root of unity. Then the position of $ \mu^i $ in the 
alcove geometry with respect to $ l $ is equal to the position of 
$r  \mu^i + (r-1)\rho $ with respect to the alcove geometry defined by
$ rl $. Indeed letting $ \omega_i $ denote the $ i$'th fundamental 
weight we have that 
 
$$ \langle r \lambda + (r-1)\rho + \rho, \omega_i \rangle =  
r \langle \lambda +  \rho, \omega_i \rangle  $$
and hence the end residues of $ r \lambda + (r-1)\rho $
equal the $ r $-multiples of the residues of $ \lambda $. 
Thus increasing the size of $ r $ 
gives rise to more weights inside the alcoves and 
makes it possible to ``separate'' hyperplanes 
through $ \lambda $ coming from different residues.
 
\medskip

We choose $ r $ large enough for 
all of the operators $ \Theta $ on the weights $ \mu^i $ 
to be on the form $ \Theta^* $ or $ \Theta_* $, that is involving no
composites of such. 
Furthermore, we obtain by choosing $r$ big enough that all the occurring weights 
$ \mu^i $ are $l$-regular partitions.

\medskip

Write $ \mu^i = \sigma $ and $ \mu^{i+1} = \tau $
and let us first assume that the operator that takes $ \sigma $ to 
$ \tau $ is 
of the type $ \Theta_* $. We check that it can 
realized through a sequence of $ f_i^{(k)}$'s.

\medskip

By assumption we have that $ \tau_i - \tau_{j} - i + j \equiv 0 \mod l
\Rightarrow \sigma_i  - \sigma_j - i + j  \equiv 0 \mod l $.
Let $ I_0  \subset [1,\ldots,m ] $
consist of those indices $ i $ such that 
$ \tau_i  - i  \equiv 0 \mod l $
and let similarly $ J_0 $ be those indices $ i  $ such that
$ \sigma_i  - i \equiv 0 \mod l $.
Then $ I_0 $ determines the hyperplanes
passing through $ \tau $ that come from the residue $ 0 \mod l $ 
and similarly $ J_0 $. Thus by assumption we have $ I_0 \subset J_0 $.

\medskip

Let $ j \in I_0 \setminus J_0 $ be minimal. We can now 
add nodes to the $j$'th line of $ \sigma $ until the last residue becomes 
$ 0 $. Each node added on the way does not give rise to 
any new coinciding residues, since otherwise there would 
be a hyperplane separating $ \sigma $ and 
$ \tau $. Similarly we deal with the other  
elements of $  \in I_0 \setminus J_0 $. But adding such nodes 
corresponds to the operation of $ f_i $ where $ i $ is 
the residue of the node. 
The other elements of 
$ j \in I_0 \setminus J_0 $ are dealt with similarly.

\medskip

At this stage, $ \sigma $ and $ \tau $ are in the same facette and 
we can add or subtract nodes to $ \sigma $, without producing coinciding 
residues, to arrive at $ \tau $. Adding 
these nodes corresponds to the operation of certain $ f_i$'s while subtracting of 
nodes corresponds to certain $ e_i $'s. But using the relations of 
$ U_q( \widehat{ {\mathfrak  sl }_l } ) $, these cancel out and
we are done in this case.

\medskip

Assume now that the operator that takes $ \mu^{i} $ 
to $ \mu^{i+1} $
is of the type $ \Theta^* $ and 
write $ \mu^i = \tau $ and $ \mu^{i+1} = \sigma $. 
Thus, there is a root $ \alpha $ such that 
$ \langle \tau + \rho, \alpha \rangle \equiv 0 $ while 
$ \langle \sigma + \rho, \alpha \rangle \not\equiv 0 $.
We may assume that $ \alpha $ is the 
only such root, by otherwise passing to a larger $l$.
Write $ \alpha = \omega_k - \omega_l $ for 
$ k < l $ and assume wlog. that the residues of 
the $ k $'th and $l$'th line of $\tau $ are $0$. 
Let $ I_0$ be as before, i.e. $I_0$ defines the 
hyperplanes passing through $ \tau $ coming from 
the residue $ 0$. Then the end residues of the lines
in $\sigma $ of indices $I_0^1 := \{i  \in I_0 | \, i \leq k \, \} $ are constant and so 
are the end residues of the lines in $ \sigma $ of indices      
$I_0^2 :=   \{i  \in I_0 |\,  i \geq l \, \} $. Let the first constant be $ n_1 $ and the second
be $ n_2 $. Assume first that $n_1 = 1 $ and $n_2 = 0 $. 
Then using Lemma 6.2 of [LLT] one checks that $ f_1^{( | I_0^1 | )} $ takes $ \tau $ to
$ \sigma $. 

\medskip

For larger values of $ n_1 $ we instead operate with the composite 
$ f_{n_1}^{( | I_0^1 | ) } \cdots   f_2^{( | I_0^1 | )}   f_1^{( | I_0^1 | )} $ on $ \tau$ 
and for larger values of $ n_2 $ we first operates with 
$ f_{n_2}^{( | I_0 | ) } \cdots   f_{2}^{( | I_0 | )}   f_1^{( | I_0 | )} $ 
on $\tau $
and then with a sequence of the first type. 
    
\medskip

This finishes the proof of the Theorem.

\begin{flushright} $ \Box $ \end{flushright}

\medskip

{\bf Remark.} The paper [GW1] by Goodman and Wenzl contains a path algorithm for affine 
Kazhdan-Lusztig polynomials valid for all Lie types. This gives a different proof of the Theorem.

\medskip

Using the argument from the beginning of section 5 
we may now conclude that the LLT-algorithm calculates $[ T(\lambda), \Delta(\mu) ] $
for $ \lambda $ a regular partition. 
But by a $q$-analogue of Karin Erdmann's result in [E], this number is equal to the decomposition number 
$d_{\lambda \mu } $ for Hecke algebras at an $ l$'th root of
unity. In other words, as claimed in the introduction of our paper, the LLT-conjecture follows from Soergel's 
algorithm.

\medskip

{\bf Remark.} The LLT conjecture only treats canonical basis coefficients $ d_{\lambda, \mu}(q) $ for
$ \lambda $ an $ l$-regular partition. On the other hand, the singular combinatorics defined in the present paper should 
work for arbitrary $ \lambda $ as well and produce decomposition numbers for the $q$-Schur algebra. In [LT1], canonical 
basis coefficients 
$ d_{\lambda, \mu}(q) $ were defined for arbitrary $\lambda$ and it was conjectured that their values at $1$
coincide with these decomposition numbers. This conjecure was proved in [VV].

\end{document}